\documentclass[conference]{IEEEtran}
\IEEEoverridecommandlockouts
\usepackage{cite}
\usepackage{url}
\usepackage{hyperref}

\usepackage{amsmath,amssymb,amsfonts}
\usepackage{algorithmic}
\usepackage{graphicx}
\usepackage{textcomp}
\usepackage{xcolor}
\usepackage{booktabs}
\usepackage{float}

\def\BibTeX{{\rm B\kern-.05em{\sc i\kern-.025em b}\kern-.08em
    T\kern-.1667em\lower.7ex\hbox{E}\kern-.125emX}}

\begin{document}

\title{Optimal Deployment of Electric Aircraft for Canadian Domestic Flights\\

\thanks{This research was funded by the Natural Sciences and Engineering Research Council of Canada (NSERC) through Alliance Mission Grants.}
}

\author{\IEEEauthorblockN{1\textsuperscript{st} Elham Soufiani}
\IEEEauthorblockA{\textit{dept. of Management Science and Engineering} \\
\textit{University of Waterloo}\\
Waterloo, ON, Canada \\
Elham.Soufiani@uwaterloo.ca}
\and
\IEEEauthorblockN{2\textsuperscript{nd} Mehrdad Pirnia}
\IEEEauthorblockA{\textit{dept. of Management Science and Engineering} \\
\textit{University of Waterloo}\\
Waterloo, ON, Canada \\
mpirnia@uwaterloo.ca}
}

\maketitle

\begin{abstract}
This paper presents a multi-period mixed-integer linear programming (MILP) framework for planning the transition from conventional to electric aircraft in regional aviation. The model jointly optimizes fleet acquisition, infrastructure deployment, and service allocation over time, while accounting for policy constraints such as emissions reduction targets, electric service share, and budget limits. A real-world case study based on Helijet’s short-haul network in British Columbia demonstrates the applicability of the model. The results show that electrification can reduce emissions by more than 70\% within five years while remaining economically viable. However, the transition is primarily limited by the capacity of the fleet and operational structure, rather than the charging infrastructure, leading to unmet demand under direct aircraft replacement. These findings emphasize the need for coordinated planning across fleet sizing, scheduling, and route prioritization to ensure a practical and efficient transition to electric aviation.

\end{abstract}


\section{INTRODUCTION}

\IEEEPARstart{T}{HE} global aviation industry currently contributes approximately 2.5\% of worldwide carbon emissions \cite{ourworldindata2024}. With international traffic projected to grow at a compound annual growth rate of 3.4\%, reaching 17.7 billion passengers over the next two decades and domestic demand expected to grow by 2.8\% by 2043 \cite{aci2025}, this environmental impact is likely to increase significantly if left unaddressed. In response, the Government of Canada has committed to achieving net-zero greenhouse gas emissions by 2050, as outlined in its Aviation Climate Action Plan (2022--2030) \cite{tc2022}. A key component of this transition is the adoption of sustainable aviation technologies, particularly electric aircraft (EA), which are most practical for short-range regional operations.

Despite this promise, technical limitations remain a major barrier. Current battery technologies have significantly lower energy density than conventional aviation fuel (approximately 250~Wh/kg versus 12,000~Wh/kg) \cite{barzkar2022,dornGomba2020}, which restricts the feasibility of electric long-haul flights. Regional aviation, however, presents a more immediate opportunity, especially in Canada, where more than half of domestic flights are less than one hour. Even within this range, transitioning to a regional fleet is far from straightforward. It involves a combination of operational, infrastructural, and financial challenges that must be addressed jointly. Battery degradation reduces usable energy over time and increases replacement costs, while charging duration and turnaround requirements limit scheduling flexibility and airport capacity \cite{airportmicrogrid2022}. These factors are closely interconnected, making fleet electrification a multi-dimensional planning problem rather than a purely technical one.

Recent research on electric aviation has largely focused on isolated subproblems. The design of the charging infrastructure has been studied using network optimization to locate and size charging stations \cite{milp,trc}, and short-term scheduling models have addressed daily flight assignment and energy dispatch for EA \cite{10345712}. Separately, battery degradation has been modeled to estimate its effect on aircraft range and lifespan under varying operational conditions \cite{lifetime}. Although these contributions advance individual aspects of electrification, they treat fleet composition, infrastructure availability, and service allocation as fixed inputs rather than joint decision variables. Strategic questions, such as when to retire conventional aircraft, which airports to electrify first, how many chargers to install, and which routes to electrify in which year, remain largely unaddressed in an integrated, multi-period framework. Without coordinating these decisions over time, operators risk either under-investing in infrastructure relative to fleet growth or retiring conventional capacity faster than EA can replace it, both of which undermine the feasibility of policy-driven transition targets.

To address this gap, this paper proposes a multi-period mixed-integer linear programming (MILP) framework that jointly optimizes fleet acquisition, airport electrification, charger deployment, and service allocation over a multi-year horizon. Unlike existing approaches, the model integrates investment and operational decisions under policy-driven constraints, including minimum electric passenger share, emissions reduction targets, annual capital budgets, and limits on fleet utilization and charging capacity. The model leverages fleet schedules and route structures. Battery degradation is incorporated as an endogenous cost component, linking operational intensity with long-term economic performance. The framework is applied to Helijet International’s harbour-to-harbour network in British Columbia, using the Beta Technologies ALIA A250 eVTOL as the candidate replacement platform. The model produces a year-by-year transition plan and shows that emission reductions of up to 72\% are achievable within the planning horizon. More importantly, the results reveal two key structural insights: first, the capacity gap between electric and conventional aircraft limits the pace of fleet transition without creating unmet demand; and second, prioritizing high-frequency routes is essential for achieving cost-effective electrification.

\section{METHODOLOGY}
The proposed planning model determines how fleet electrification, airport electrification, and charging infrastructure evolve over the planning horizon $T$, while optimally allocating electric and conventional services. The model considers planning years $t \in T$, itineraries $i \in I$, airports $l \in L$, and paired itineraries $q=(i_d(q), i_r(q)) \in Q$. The parameters and decision variables are summarized in Table~\ref{tab:notation}.

\begin{table}[!htbp]
\footnotesize
\caption{Compact notation for the multi-period fleet transition MILP.}
\label{tab:notation}
\centering
\footnotesize
\renewcommand{\arraystretch}{1.05}
\begin{tabular}{p{0.19\columnwidth} p{0.77\columnwidth}}
\hline
\multicolumn{2}{c}{\textbf{Parameters}} \\
\hline
$D_{i,t}$ & Passenger demand for itinerary $i$ in year $t$ (passengers/year). \\
$F_{i,t}$ & Maximum number of flights available on itinerary $i$ in year $t$ (flights/year). \\
$K^E,K^C$ & Seats per electric and conventional aircraft, respectively (seats/aircraft). \\
$\tau_i^E,\tau_i^C$ & Flight-hours required by one electric/conventional flight on itinerary $i$ (hr/flight). \\
$\bar{\tau}_t^E,\bar{\tau}_t^C$ & Annual available operating hours per electric/conventional aircraft (hr/aircraft/year). \\
$e_{i,t}$ & Charging energy required for one electric flight on itinerary $i$ in year $t$ (kWh/flight). \\
$\eta^{ch}$ & Charging efficiency (dimensionless, $0 < \eta^{ch} \leq 1$). \\
$\bar{e}_{q,l,t}$ & Charging energy attributed to airport $l$ for one electric paired operation of $q$ in year $t$ (kWh/operation). \\
$P_{l,t}$ & Annual charging throughput provided by one charger at airport $l$ in year $t$ (kWh/charger/year). \\
$\bar{m}_l$ & Maximum number of chargers allowed at airport $l$ (chargers). \\
$\bar{h}_{l,t}$ & Maximum number of chargers that can be added at airport $l$ in year $t$ (chargers/year). \\
$\Phi_{q,t}^1$ & Equals $1$ if route pair $q$ is feasible with at least one electrified endpoint; 0 otherwise (binary). \\
$\Phi_{q,t}^2$ & Equals $1$ if route pair $q$ is feasible only when both endpoints are electrified; 0 otherwise (binary). \\
$\alpha_t$ & Minimum required electric service share in year $t$ (dimensionless, $0 \leq \alpha_t \leq 1$). \\
$\rho_t$ & Minimum required emissions reduction fraction in year $t$ (dimensionless, $0 \leq \rho_t \leq 1$). \\
$EM_t^{base}$ & Baseline emissions in year $t$ (kgCO$_2$e/year). \\
$E_{i,t}^E,E_{i,t}^C$ & Emissions per electric/conventional flight on itinerary $i$ in year $t$ (kgCO$_2$e/flight). \\
$BUD_t$ & Investment budget in year $t$ (\$/year). \\
$C_t^{buy,E}$ & Purchase cost of one EA in year $t$ (\$/aircraft). \\
$C_{l,t}^{hub}$ & Cost of electrifying airport $l$ in year $t$ (\$/airport). \\
$C_{l,t}^{ch}$ & Cost of installing one charger at airport $l$ in year $t$ (\$/charger). \\
$C_t^{fix,E},C_t^{fix,C}$ & Annual fixed cost per electric/conventional aircraft (\$/aircraft/year). \\
$c_t^{hr,E},c_t^{hr,C}$ & Operating cost per electric/conventional flight-hour (\$/hr). \\
$c_{l,t}^{grid}$ & Grid-electricity cost per unit of charging energy at airport $l$ in year $t$ (\$/kWh). \\
$\gamma_t$ & Battery degradation cost per unit of electric energy usage in year $t$ (\$/kWh). \\
$C_{i,t}^{pen}$ & Penalty cost for one unit of unmet demand on itinerary $i$ in year $t$ (\$/passenger). \\
\hline
\multicolumn{2}{c}{\textbf{Decision variables}} \\
\hline
$b_t^E$ & Number of EA purchased in year $t$ (aircraft). \\
$r_t^C$ & Number of conventional aircraft retired in year $t$ (aircraft). \\
$N_t^E,N_t^C$ & Electric and conventional fleet sizes in year $t$ (aircraft). \\
$z_{l,t}$ & 1 if airport $l$ is electrified in year $t$; 0 otherwise (binary). \\
$w_{l,t}$ & 1 if airport $l$ is electrified by the end of year $t$; 0 otherwise (binary). \\
$h_{l,t}$ & Number of chargers installed at airport $l$ in year $t$ (chargers). \\
$m_{l,t}$ & Number of chargers available at airport $l$ by the end of year $t$ (chargers). \\
$y_{q,l,t}$ & 1 if airport $l$ is selected as a charging-support airport for pair $q$ in year $t$; 0 otherwise (binary). \\
$y_{q,t}^1$ & 1 if paired itinerary $q$ is activated under the single-end charging feasibility mode in year $t$; 0 otherwise (binary). \\
$y_{q,t}^2$ & 1 if paired itinerary $q$ is activated under the dual-end charging feasibility mode in year $t$; 0 otherwise (binary). \\
$y_{q,t}$ & 1 if paired itinerary $q$ is allowed to operate electrically in year $t$; 0 otherwise (binary). \\
$v_{q,l,t}$ & Number of electric paired operations of $q$ assigned to airport $l$ in year $t$ (operations). \\
$x_{i,t}^E,x_{i,t}^C$ & Number of electric/conventional flights operated on itinerary $i$ in year $t$ (flights/year). \\
$n_{i,t}^E,n_{i,t}^C$ & Demand served electrically/conventionally on itinerary $i$ in year $t$ (passengers/year). \\
$u_{i,t}$ & Unmet demand on itinerary $i$ in year $t$ (passengers/year). \\
\hline
\end{tabular}
\end{table}

\subsection{Objective function}

\vspace{-0.6cm}
\begin{align}
\min \; Z = \sum_{t\in T} \Bigg[
& C_t^{buy,E} b_t^E
+ \sum_{l\in L} C_{l,t}^{hub} z_{l,t}
+ \sum_{l\in L} C_{l,t}^{ch} h_{l,t}\nonumber\\
&+ C_t^{fix,E} N_t^E
+ C_t^{fix,C} N_t^C \nonumber\\
&+ c_t^{hr,E} \sum_{i\in I} \tau_i^E x_{i,t}^E
+ c_t^{hr,C} \sum_{i\in I} \tau_i^C x_{i,t}^C \nonumber\\
&+ \sum_{l\in L} c_{l,t}^{grid}
\sum_{q\in Q}\frac{\bar e_{q,l,t}}{\eta^{ch}}\, v_{q,l,t}
+ \gamma_t \sum_{i\in I} e_{i,t} x_{i,t}^E \nonumber\\
&+ \sum_{i\in I} C_{i,t}^{pen} u_{i,t}
\Bigg].
\end{align}

\noindent
The objective minimizes the total cost of the aviation system over the planning horizon. It includes capital investments in EA, $C_t^{buy,E}$, airport electrification, $C_{l,t}^{hub}$, and charging infrastructure, $C_{l,t}^{ch}$, costs, along with fixed fleet costs, $C_t^{fix,E}, C_t^{fix,C}$. Operating costs are modeled according to flight hours for electric and conventional services, $x_{i,t}^E, x_{i,t}^C$, with the corresponding durations $\tau_i^E$ and $\tau_i^C$. 

Electricity costs are captured through charging energy consumption, adjusted by charging efficiency, while battery degradation is represented through the term $\gamma_t$. Finally, unmet demand is penalized through $u_{i,t}$ to ensure service reliability.


\subsection{Fleet evolution}

\vspace{-0.5cm}
\begin{align}
N_t^E &= N_{t-1}^E + b_t^E, && \forall t, \\
N_t^C &= N_{t-1}^C - r_t^C, && \forall t, \\
0 & \le r_t^C \le N_{t-1}^C, && \forall t.
\end{align}

\noindent
These constraints describe the evolution of electric, $N_t^E$, and conventional, $N_t^C$, fleet sizes over time. EA are added through annual purchases $b_t^E$, while conventional aircraft are gradually retired through $r_t^C$. The constraint $r_t^C \le N_{t-1}^C$ ensures that retirements do not exceed the existing fleet size.


\subsection{Airport electrification and charger expansion}

\vspace{-0.5cm}
\begin{align}
w_{l,t} &= w_{l,t-1} + z_{l,t}, && \forall l,t, \\
z_{l,t} &\le 1 - w_{l,t-1}, && \forall l,t, \\
m_{l,t} &= m_{l,t-1} + h_{l,t}, && \forall l,t, \\
m_{l,t} &\le \bar{m}_l\, w_{l,t}, && \forall l,t, \\
0 & \le h_{l,t} \le \bar{h}_{l,t}, && \forall l,t.
\end{align}

\noindent
These constraints model the evolution of airport electrification and charging infrastructure over time. The binary variable $z_{l,t}$ determines whether the airport $l$ is electrified in year $t$, while $w_{l,t}$ monitors its cumulative electrification status. 

The charger deployment is represented by $h_{l,t}$ (new installations) and $m_{l,t}$ (total installed chargers). Chargers can only be installed at electrified airports, enforced by $m_{l,t} \le \bar{m}_l w_{l,t}$, and annual installations are limited by $\bar{h}_{l,t}$.


\subsection{Route-type electric feasibility}

\vspace{-0.5cm}
\begin{align}
\Phi_{q,t}^{1} + \Phi_{q,t}^{2} &\le 1,
&& \forall q, t \label{eq:phi_type}
\\
y_{q,t}^{1} &\le \Phi_{q,t}^{1},
&& \forall q,t, \label{eq:y1_phi}
\\
y_{q,t}^{1} &\le w_{o(q),t} + w_{d(q),t},
&& \forall q,t, \label{eq:y1_upper}
\\
2y_{q,t}^{1} &\ge \Phi_{q,t}^{1}\big(w_{o(q),t}+w_{d(q),t}\big),
&& \forall q,t, \label{eq:y1_lower}
\\
y_{q,t}^{2} &\le \Phi_{q,t}^{2},
&& \forall q,t, \label{eq:y2_phi}
\\
y_{q,t}^{2} &\le w_{o(q),t},
&& \forall q,t, \label{eq:y2_upper_o}
\\
y_{q,t}^{2} &\le w_{d(q),t},
&& \forall q,t, \label{eq:y2_upper_d}
\\
y_{q,t}^{2} &\ge \Phi_{q,t}^{2}\big(w_{o(q),t}+w_{d(q),t}-1\big),
&& \forall q,t, \label{eq:y2_lower}
\\
y_{q,t} &= y_{q,t}^{1} + y_{q,t}^{2},
&& \forall q,t. \label{eq:y_total}
\end{align}

\noindent
These constraints determine whether a paired itinerary $q$ (outbound and return flights) can be operated electrically based on both technical feasibility and infrastructure availability. The binary parameters $\Phi_{q,t}^{1}$ and $\Phi_{q,t}^{2}$ classify each paired route according to their charging requirements. Specifically, $\Phi_{q,t}^{1}=1$ indicates that the route can be operated with at least one electrified endpoint, while $\Phi_{q,t}^{2}=1$ requires both endpoints to be electrified.

The variables $y_{q,t}^{1}$ and $y_{q,t}^{2}$ activate the corresponding feasibility mode based on the electrification status of the origin and destination airports, i.e., $w_{o(q),t}$ and $w_{d(q),t}$. The variable $y_{q,t}$ indicates whether the paired route $q$ is eligible for electric operation in year $t$.

This formulation ensures that electric service is allowed only when the required charging configuration, either single-end or dual-end, is satisfied, without imposing unnecessary infrastructure requirements.

\subsection{Charging-support assignment on paired itineraries}

\vspace{-0.5cm}
\begin{align}
y_{q,l,t} &= 0,
&& \forall q,t,\\
& && l\notin\{o(q),d(q)\},\nonumber \\
y_{q,l,t} &\le w_{l,t},
&& \forall q,l,t,\\
\sum_{l\in L} y_{q,l,t}
&= y_{q,t}^1 + 2y_{q,t}^2,
&& \forall q,t,\\
v_{q,l,t} &\le F_{i_d(q),t}\, y_{q,l,t},
&& \forall q,l,t,\\
\sum_{l\in L} v_{q,l,t}
&= (\Phi_{q,t}^1 + 2\Phi_{q,t}^2)\, x_{i_d(q),t}^E,
&& \forall q,t.
\end{align}

\noindent
The constraints in this block determine how charging support is allocated between the endpoints of each paired itinerary. Charging can only be provided at the airports of origin or destination, only if those airports are electrified.

Depending on the route type, either one or both endpoints must provide charging support. The binary variables $y_{q,l,t}$ identify the endpoints selected to provide this support. The integer variable $v_{q,l,t}$ then allocates the number of electric paired operations supported at each airport, ensuring consistency between electric operations and available charging infrastructure.


\subsection{Electric-flight consistency on paired itineraries}

\vspace{-0.5cm}
\begin{align}
x_{i_d(q),t}^E &\le F_{i_d(q),t}\, y_{q,t}, && \forall q,t,\\
x_{i_r(q),t}^E &\le F_{i_r(q),t}\, y_{q,t}, && \forall q,t, \\
x_{i_d(q),t}^E &= x_{i_r(q),t}^E && \forall q,t.
\end{align}

\noindent
These constraints ensure consistency between the outbound and return legs of each paired itinerary under electric operation. When a route is enabled for electrification, $y_{q,t}=1$, both legs must operate the same number of electric flights, ensuring that the aircraft can complete round-trip flights. The number of electric flights on each leg is further limited by the available scheduled capacity $F_{i,t}$.


\subsection{Demand satisfaction and flight capacity}

\vspace{-0.5cm}
\begin{align}
n_{i,t}^E + n_{i,t}^C + u_{i,t} &= D_{i,t}, && \forall i,t, \\
n_{i,t}^E &\le K^E x_{i,t}^E, && \forall i,t, \\
n_{i,t}^C &\le K^C x_{i,t}^C, && \forall i,t, \\
x_{i,t}^E + x_{i,t}^C &= F_{i,t} && \forall i,t.
\end{align}

\noindent
Here, the constraints ensure that passenger demand is properly accounted for on each route. Demand $D_{i,t}$ is served by electric flights, $n_{i,t}^E$, conventional flights, $n_{i,t}^C$, or remains unmet, $u_{i,t}$.

Seat-capacity constraints link the number of passengers with the number of flights through aircraft capacities $K^E$ and $K^C$. Finally, the total number of flights assigned to each route is limited by the available schedule $F_{i,t}$.


\subsection{Fleet operating-hour limits}

\vspace{-0.5cm}
\begin{align}
\sum_{i} \tau_i^E x_{i,t}^E &\le \bar{\tau}_t^E N_t^E,  && \forall t,\\
\sum_{i} \tau_i^C x_{i,t}^C &\le \bar{\tau}_t^C N_t^C  && \forall t.
\end{align}

\noindent
These constraints limit total operations according to the available fleet capacity. The total flight hours required, given by the number of flights $x_{i,t}$ and their durations $\tau_i$, must not exceed the annual operating time available from the fleet, represented by $\bar{\tau}_t^E N_t^E$ and $\bar{\tau}_t^C N_t^C$ for electric and conventional aircraft, respectively.


\subsection{Charging-capacity limits}

\vspace{-0.5cm}
\begin{align}
\sum_{q\in Q} \frac{\bar e_{q,l,t}}{\eta^{ch}}\, v_{q,l,t}
\le P_{l,t} m_{l,t},
&& \forall l,t.
\end{align}

\noindent
This constraint ensures that the total charging demand at each airport does not exceed the available charging capacity. The left-hand side represents the energy required to support electric operations assigned to airport $l$, adjusted for charging efficiency, while the right-hand side captures the total capacity provided by the installed chargers.

This formulation directly links operational decisions with infrastructure availability, ensuring that electric flight schedules remain feasible given the charging resources at each airport.


\subsection{Policy and budget constraints}

\vspace{-0.5cm}
\begin{align}
&\sum_{i} n_{i,t}^E 
\ge \alpha_t \sum_{i} D_{i,t},&& \forall t, \\
&\sum_{i} E_{i,t}^E x_{i,t}^E
+ \sum_{i} E_{i,t}^C x_{i,t}^C
\le (1-\rho_t) EM_t^{base} ,   && \forall t,\\
& C_t^{buy,E} b_t^E
+ \sum_{l} C_{l,t}^{hub} z_{l,t}
+ \sum_{l} C_{l,t}^{ch} h_{l,t}
\le BUD_t  && \forall t.
\end{align}

\noindent
This block imposes policy and financial requirements on the system. The first constraint enforces a minimum share of the demand to be served by electric operations, defined by $\alpha_t$. The second ensures that total emissions from both electric and conventional flights remain below a target level relative to baseline emissions $EM_t^{base}$. 

The third constraint limits the total annual investment in fleet acquisition and infrastructure expansion to the available budget $BUD_t$.


\subsection{Variable domains}

\vspace{-0.5cm}
\begin{align}
&b_t^E,r_t^C,N_t^E,N_t^C,h_{l,t},m_{l,t},x_{i,t}^E,x_{i,t}^C,v_{q,l,t}\in\mathbb{Z}_+,\\
&n_{i,t}^E,n_{i,t}^C,u_{i,t}\ge 0,\\
&w_{l,t},z_{l,t},y_{q,l,t},y_{q,t}^1,y_{q,t}^2,y_{q,t}\in\{0,1\}.
\end{align}

\noindent
The above constraints define the domains of the decision variables.

\section{CASE STUDY OVERVIEW AND RESULTS}


To demonstrate the applicability of the proposed MILP framework, the model is applied to Helijet International, Canada’s only scheduled passenger helicopter airline. Its compact and high-frequency network connecting Vancouver, Victoria, and Nanaimo using Sikorsky S-76 helicopters, characterized by short legs and fixed harbour-to-harbour routes, provides a suitable setting for evaluating EA integration.

The current fleet consists of S-76 twin-turbine helicopters with 12 seats and an emission rate of 568~kgCO$_2$e per flight-hour~\cite{helijet2025schedule}. The electrification scenario considers the Beta Technologies ALIA A250 eVTOL, a five-passenger aircraft with a 460~km range~\cite{beta_alia_specs}. It is powered by a 225~kWh battery and recharged in under one hour using a 320~kW DC fast charger (such as Charge Cube)~\cite{beta2024network}. A baseline energy intensity of 0.650~kWh/km is assumed, improving by 1\% annually. The A250 is selected due to its vertical takeoff capability required by Helijet’s harbour infrastructure.

\begin{figure}[H]
    \centering
    \includegraphics[width=.7\linewidth]{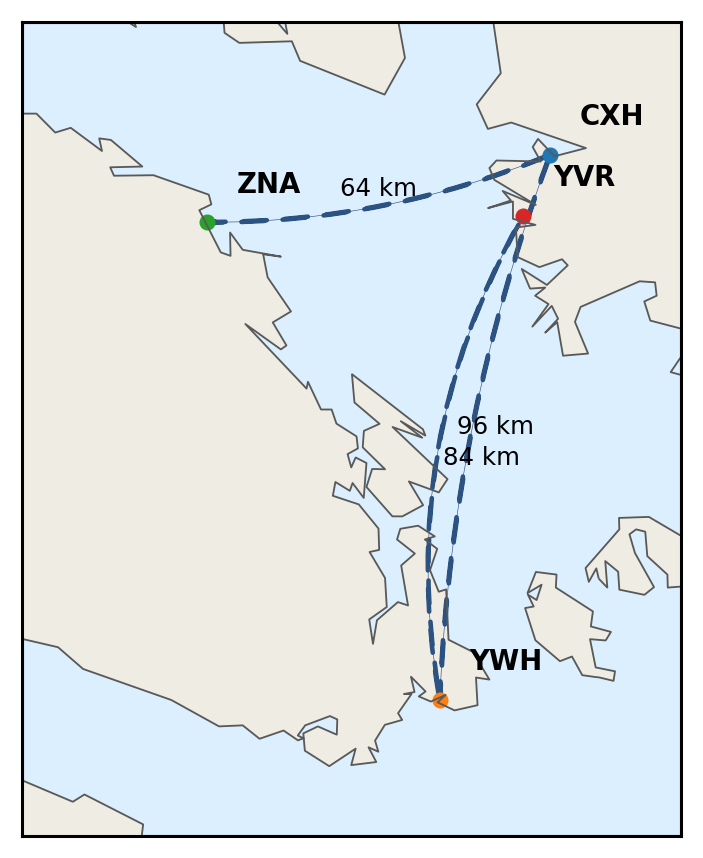}
    \caption{Helijet route map}
    \label{fig:heli_map}
\end{figure}

Three paired itineraries are modeled according to Helijet’s 2024--2026 schedules~\cite{helijet2025schedule}. Route distances are computed between heliports, with conventional flight times from published timetables and electric times estimated from A250 performance. Eligibility is confirmed using the Fall/Winter 2024--25 schedule~\cite{helijet2024schedule}. The Vancouver Harbour--Victoria Harbour route is the primary corridor, while the Nanaimo route serves a secondary market, as shown in Fig.~\ref{fig:heli_map}. The Vancouver Airport--Victoria route is included as a reference case with reduced and declining demand following its direct flight suspension. All routes satisfy single-endpoint charging feasibility, as round-trip energy remains within usable battery limits.

The baseline parameters for the initial year (2026) and their evolution over time are selected to reflect both current operational conditions and anticipated technology and policy trends. Table~\ref{tab:rates} summarizes the key rates used in the model. EA characteristics, including range, battery capacity (225~kWh), and charging capability, are based on publicly available specifications of the ALIA A250 platform, which is designed for short-haul regional operations and fast-turnaround charging~\cite{motorwatt, beta_alia_specs, beta2024network}. The assumed energy intensity of 0.650~kWh/km represents a cruise-based estimate and is treated conservatively with a gradual improvement rate.

Battery degradation is incorporated through the parameter $\gamma = 1.852~\$/\text{kWh}$, derived from battery replacement cost, usable capacity, and cycle life. This links operational intensity to long-term capital expenditure and reflects the economic impact of battery wear. Demand growth is set to 3.5\% annually, consistent with regional aviation trends~\cite{aci2025}, while policy targets for electric service share ($\alpha_t$) and emissions reduction ($\rho_t$) are introduced progressively to reflect realistic adoption pathways. Cost-related parameters, including aircraft price and grid electricity rates, are modeled with moderate annual escalation or decline to capture expected market evolution. The initial acquisition cost of the eVTOL is set to approximately \$5.5M CAD, based on publicly reported estimates of around \$4M USD for the ALIA A250 platform, adjusted to Canadian dollars~\cite{motorwatt}.

\begin{table}[!t]
\footnotesize
\caption{Key annual rates and initial assumptions}
\label{tab:rates}
\centering
\renewcommand{\arraystretch}{1.1}
\begin{tabular}{p{3.6cm} p{3.9cm}}
\toprule
\textbf{Parameter} & \textbf{Value / Trend} \\
\midrule
Demand growth & $3.5\%$ per year \\
Energy intensity & $0.650$ kWh/km ($-1\%$ per year) \\
EA acquisition cost & Decreasing ($\approx 2\%$ per year) \\
Grid electricity price & Increasing ($\approx 2\text{--}3.75\%$ per year) \\
Battery degradation ($\gamma$) & $1.852\ \$/\text{kWh}$ \\
Electric share target ($\alpha_t$) & $0\%$ to $20\%$ (ramp) \\
Emissions reduction target ($\rho_t$) & $0\%$ to $20\%$ (ramp) \\
\bottomrule
\end{tabular}
\end{table}

\begin{figure}[htbp]
\centering
\includegraphics[width=\columnwidth]{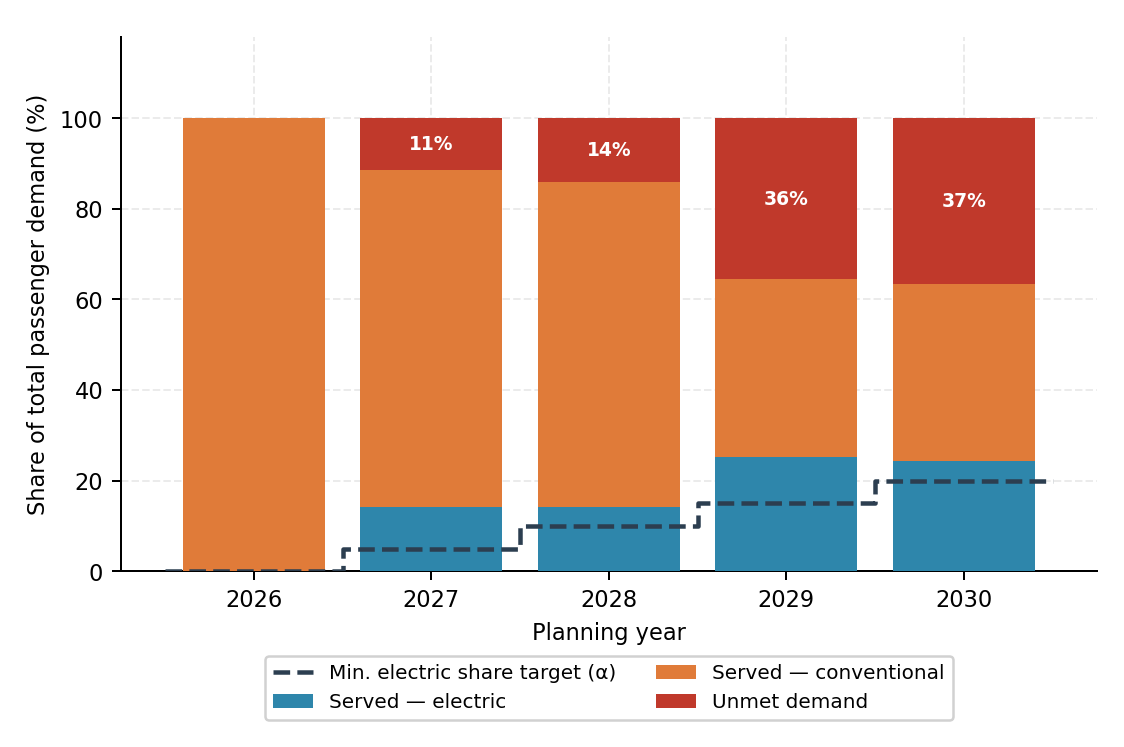}
\caption{Passenger demand allocation and electric share relative to $\alpha$.}
\label{fig:demand}
\end{figure}
Figure~\ref{fig:demand} reveals the central operational outcome. Although electric service consistently exceeds the policy target $\alpha$, the unmet demand emerging from 2027 increases to approximately 37\% by 2030. This is driven by the seat capacity gap between aircraft types (5 vs.\ 12 passengers), combined with fixed flight schedules and limited fleet expansion in the early years. This limitation is widely recognized as a key barrier to rapid electrification in regional aviation, expected to improve as manufacturers scale up capacity~\cite{icct, alia_up}. As conventional aircraft are retired, the total capacity of the system declines, and the model accepts unmet demand as a cost-minimizing trade-off.

\begin{figure}[htbp]
\centering
\includegraphics[width=\columnwidth]{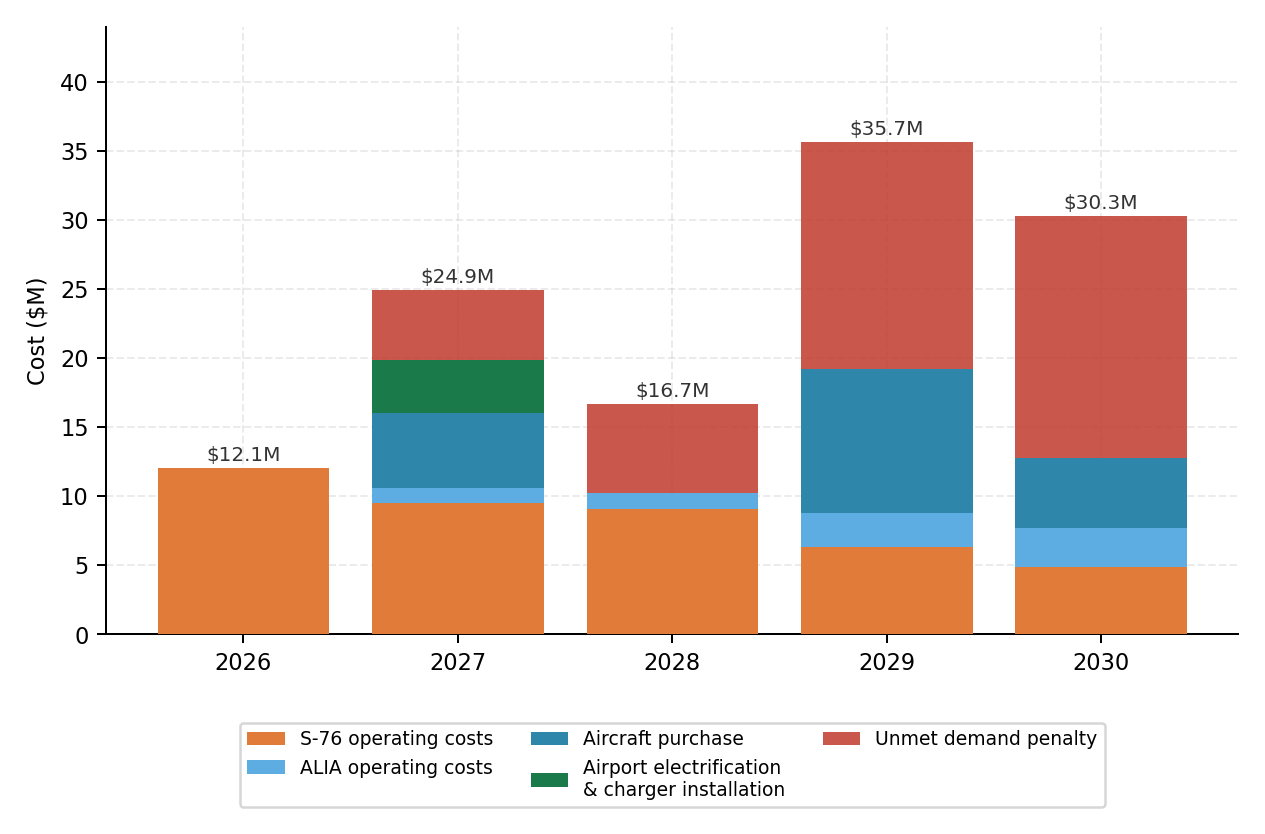}
\caption{Annual system cost breakdown.}
\label{fig:cost}
\end{figure}
The cost structure, Figure~\ref{fig:cost}, reflects this transition dynamic. Conventional operating costs dominate the early years, while 2027 represents the primary investment year due to aircraft acquisition and infrastructure deployment. In later years, unmet demand penalties become the dominant cost component, reaching more than \$17M annually. This indicates that system inefficiency is driven primarily by insufficient transport capacity rather than capital or operating costs.

The charging infrastructure is not a limiting factor. The model determines that a single high-power charger per electrified hub is sufficient to support all electric operations throughout the horizon.

\begin{figure}[htbp]
\centering
\includegraphics[width=\columnwidth]{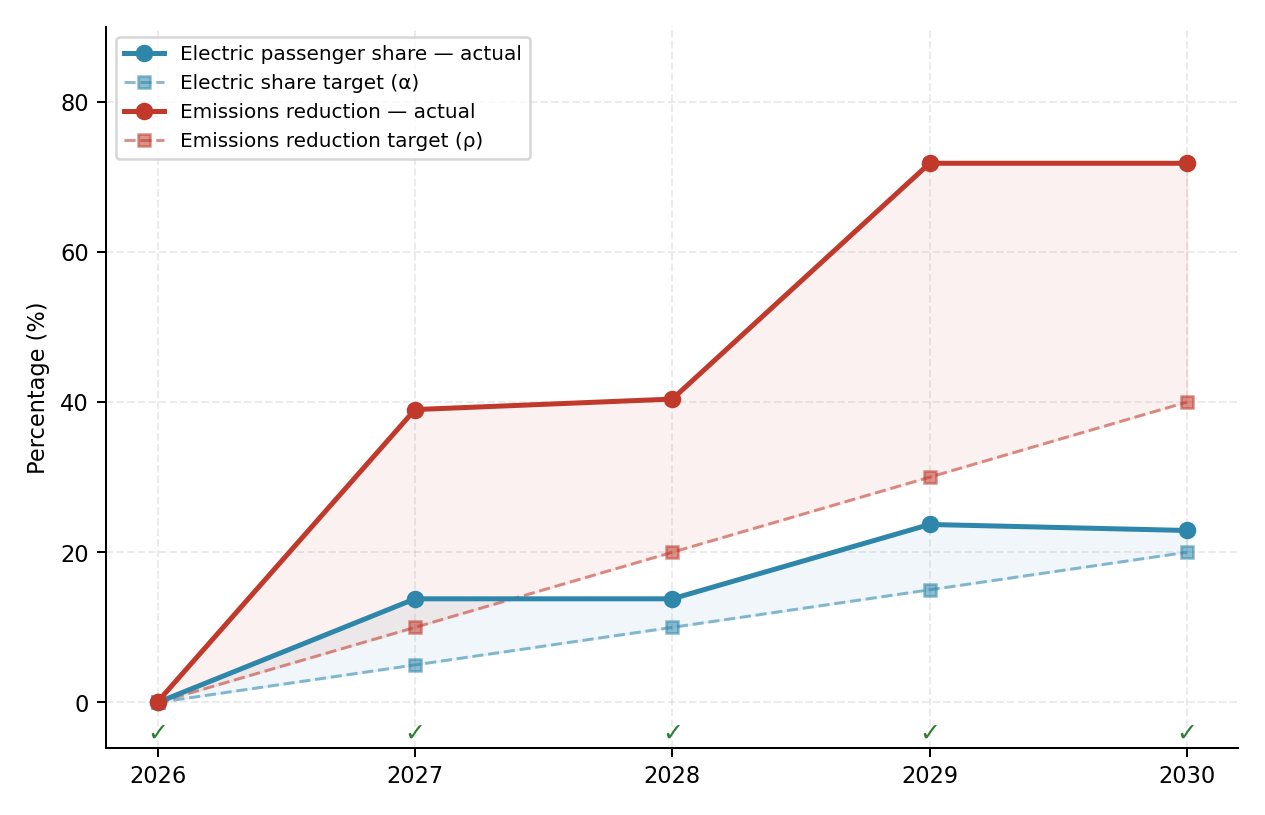}
\caption{Policy constraint performance: electric share ($\alpha$) and emissions reduction ($\rho$).}
\label{fig:policy}
\end{figure}

Figure~\ref{fig:policy} shows that despite capacity limitations, all policy constraints are satisfied. The share of eVTOLs consistently exceeds targets and emissions reductions significantly outperform requirements, reaching more than 70\% by 2029. These reductions are achieved early in the transition due to the concentration of electric operations on the route with the highest-impact and the low carbon intensity of the electricity supply. The key decision outcomes are summarized in Table~\ref{tab:keydecisions}.

\begin{table}[!h]
\footnotesize
\caption{Key decision outcomes (non-visualized)}
\label{tab:keydecisions}
\centering
\renewcommand{\arraystretch}{1.1}
\begin{tabular}{p{3cm} p{4.5cm}}
\toprule
Decision & Outcome \\
\midrule
Aircraft purchases & 1 (2027), 2 (2029), 1 (2030) \\
Conventional retirements & 1 (2028), 1 (2029), 3 (2030) \\
Electrified airports & CXH and YWH (from 2027) \\
Chargers installed & 1 per hub; no expansion required \\
Route electrification & Only CXH$\leftrightarrow$YWH electrified \\
Budget constraint & Non-binding in all years \\
\bottomrule
\end{tabular}
\end{table}

\section{DISCUSSION}

The results indicate that electrification is primarily constrained by the fleet capacity and operational structure rather than the infrastructure or energy supply. The central challenge arises from the capacity mismatch between conventional and EA under fixed schedules. Concentrating electrification on high-frequency corridors delivers immediate cost and emission benefits, supporting a targeted rather than network-wide deployment strategy. However, the fleet transition must be carefully managed. A one-to-one replacement approach leads to unmet demand and increasing penalty costs, highlighting the need for complementary measures such as increased flight frequency, larger-capacity EA, or delayed retirement of conventional units.

Flexibility in scheduling provides an important operational lever. The shorter flight times of EA enable higher aircraft utilization, allowing more flights to be scheduled as adoption increases. In addition, improvements in battery technology, particularly higher charging power and reduced charging time, can shorten ground turnaround, further increasing feasible flight frequency and partially mitigating capacity loss without requiring additional fleet investment. In contrast, the charging infrastructure itself is not a limiting factor. The system operates well below installed capacity, indicating that additional chargers would yield limited benefit relative to improvements in fleet capacity and scheduling.

\section{CONCLUSIONS and FUTURE WORK}

This study presented a multi-period MILP framework for planning the transition from conventional to EA in regional aviation systems. The model is applied to the Helijet network, which demonstrated that electrification can deliver substantial emissions reductions, exceeding 70\% within the planning horizon, while remaining economically viable. The results showed that electrification is not limited by charging infrastructure or energy availability, but primarily by fleet capacity and operational structure. In particular, the capacity mismatch between conventional and EA on fixed schedules leads to unmet demand.

The findings suggest that effective electrification requires more than direct substitution of EA. Improvements in aircraft utilization, scheduling flexibility, and charging turnaround times can significantly enhance system performance. Future work will extend the model to incorporate dynamic scheduling decisions, stochastic demand, and evolving aircraft technologies with higher capacity and faster charging capabilities. In addition, integrating policy design with service-level constraints and exploring network-wide electrification strategies will provide further insight into scalable and operationally feasible transition pathways.

\section{ACKNOWLEDGMENT}

The authors would like to thank their colleagues, reviewers, and editors for their constructive feedback, which has significantly improved the quality of this work. The authors also gratefully acknowledge Benjamin Singer and Elika Farokhi for their valuable assistance and contributions throughout this project. This work was funded by the Natural Sciences and Engineering Research Council of Canada (NSERC) through the Alliance Mission Grant program and supported by the Waterloo Institute for Sustainable Aeronautics (WISA).

\bibliographystyle{IEEEtran}
\bibliography{References}

\end{document}